\documentclass[smallextended,a4paper]{svjour3}
\journalname{Differential Equations and Dynamical Systems}
\usepackage{graphicx}
\usepackage{amssymb}
\usepackage[colorlinks,linkcolor=blue,citecolor=blue]{hyperref}
\usepackage{anyfontsize}
\usepackage{geometry}
\usepackage{graphicx}
\smartqed
\AtBeginDocument{%
\setlength{\oddsidemargin}{\dimexpr(\paperwidth-\textwidth)/2-1in}%
\setlength{\evensidemargin}{\oddsidemargin}%
\setlength{\topmargin}{%
\dimexpr(\paperheight-\textheight)/2-\headheight-\headsep-1in}}

\newtheorem{defn}{Definition}[section]
\newtheorem{lemm}{Lemma}[section]
\newtheorem{prop}{Proposition}[section]
\newtheorem{cory}{Corollary}[section]
\newtheorem{theo}{Theorem}[section]
\newtheorem{rema}{Remark}[section]

\begin{document}
%
\title{Dynamics of locally coupled agents with next nearest neighbor
interaction}
\titlerunning{Dynamics of locally coupled agents}
\author{J.~Herbrych$^{1,2,3,*}$ \and A.~G.~Chazirakis$^{4,\dagger}$
\and N.~Christakis$^{1,\dagger}$ \and J.~J.~P.~Veerman$^{1,5,\S}$}
\authorrunning{J.~Herbrych {\it et al.}} 
\institute{$^{1}$Crete Center for Quantum Complexity and Nanotechnology,\\
Department of Physics,\\
University of Crete,\\
Heraklion, Greece\\\\
$^{2}$Department of Physics and Astronomy,\\
University of Tennessee,\\
Knoxville, TN, USA\\\\
$^{3}$Materials Science and Technology Division,\\
Oak Ridge National Laboratory,\\
Oak Ridge, TN, USA\\\\
$^{4}$Department of Mathematics and Applied Mathematics,\\
University of Crete,\\
Heraklion, Greece\\\\
$^{5}$Fariborz Maseeh Department of Mathematics and Statistics,\\
Portland State University,\\
Portland, OR, USA\\\\
$^*$\email{jherbryc@utk.edu}\\
$^\dagger$\email{thazir@gmail.com}\\
$^\ddagger$\email{nchristakis@tem.uoc.gr}\\
$^\S$\email{veerman@pdx.edu}
}
\date{Received: date / Accepted: date}
\maketitle
%
\begin{abstract}\begin{small}
We consider large but finite systems of identical agents on the line with up
to next nearest neighbor asymmetric coupling. Each agent is modelled by a
linear second order differential equation, linearly coupled to up to four of
its neighbors. The only restriction we impose is that the equations are
\emph{decentralized}. In this generality we give the conditions for stability
of these systems. For stable systems, we find the response to a change of
course by the leader. This response is at least linear in the size of the
flock. Depending on the system parameters, two types of solutions have been
found: damped oscillations and \emph{reflectionless waves}. The latter is a
novel result and a feature of systems with at least next nearest neighbor
interactions. Analytical predictions are tested in numerical simulations.
\keywords{Dynamical Systems \and Chaotic Dynamics \and Optimization and
Control \and Multi-agent Systems}
\PACS{05.45.Pq \and 07.07.Tw \and 45.30.+s \and 73.21.Ac}
\end{small}\end{abstract}

\section{Introduction}
\label{chap:intro}

Coupled second order ordinary linear differential equations, coupled
oscillators for short, play an important role in almost all areas of science
and technology (see the introduction of~\cite{Lewis2014} for a recent review).
The phenomena of coupled systems appear on all length-- and time--scales: from
synchronization of power generators in power-grid networks
\cite{Backhaus2013,Motter2013}, through the traffic control of vehicular
platoons
\cite{Bamieh2002,Cruz2007,Defoort2008,Lin2012,Barooah2013,Herman2015},
collective decision-making in biological systems
\cite{Reynolds1987,Strogatz1993,Couzin2005,Attanasi2014,Pearce2014} (e.g.,
transfer of long--range information in flocks of birds), to the atomic scale
lattice vibrations (so-called phonons), just to name few of them. The nature
of communication within such a systems crucially influences the behavior of
it. In the presence of centralized information, e.g., the knowledge of the
desired velocity by members of a flock, the performance of many of these
systems is good \cite{Bamieh2002,Lin2012,Barooah2013} in the sense that the
trajectories of the agents quickly converge to coherent (or synchronized)
motion. On the other hand, in decentralized systems convergence to coherent
motion is much less obvious, since no overall goal is observed by all agents.
In this case, the only available observations (i.e., of position and/or
velocity) are relative to the agent. The complication of the problem is even
greater if information is exchanged only locally - by agents in a neighborhood
that is small in comparison to the system size.

As an aside, we point out that there exist another class of somewhat similar
problems, also with a wide range of applications, namely the dynamics of
\emph{consensus}, see~\cite{Ren,Hegselmann}. The difference is that in our
case the agents are Newtonian (are subject to force $m\ddot x$, i.e., mass
$\times$ accelerations), whereas in consensus type equation, that is not the
case. Therefore consensus equations tend to be coupled \emph{first} order
ordinary differential equations with a very different behavior. In particular,
we do not expect to see wave-like behavior in consensus equations, whereas
they in our equations those play a prominent role. In what follows we will
restrict ourselves to coupled \emph{second} order differential equations.

It is of significant importance to develop a theory that deals with systems
where agents may interact with few nearby agents. In the case of physical
systems with symmetric interactions and no damping (such as harmonic
crystals), this theory exists and can be found in textbooks
\cite{Ashcroft1976}. It consists of imposing
\emph{periodic boundary conditions}, and then asserting that the solutions of
the periodic system behave the same way as in the system with non-trivial
boundary conditions, except near the boundary. Although we know of no formal
proof in the literature that this is correct, this method of solution has been
used for about a century with great success. Neither is it the case that we
can rely on previous studies of discretization of a second order partial
differential equation. Indeed the finite difference method applied to a wave
equation with convection will give rise (for small enough mesh) to nearly
symmetric equations \cite{Haberman,VanLoan}.

In flocks there is no reason for the interactions to be symmetric or undamped,
as is the case in the study of harmonic crystals. The equations studied here
are therefore more general. Furthermore in flocks it is desirable to have a
two parameter set of equilibria, namely motion with constant velocity and
constant distance between any two consecutive agents (\emph{coherent motion}).
Thus it is necessary to study a more general problem, namely convergence to
coherent motion in the presence of asymmetry and damping. In this paper we
generalize the previously successful approach (periodic boundary conditions).
In doing this one needs to be aware though that (i) assumptions or conjectures
needed to solve the old problem must be investigated again as they may not be
justified anymore, and (ii) new phenomena may arise. For more details
see~\cite{Cantos2014-1,Cantos2014-2,Herman2015}.

In the case of linear response theory in solid state physics
\cite{Ashcroft1976}, when a system of symmetrically coupled undamped
oscillators is perturbed, the signal will typically travel through the entire
system at constant velocity without damping. In our case, the system is
generally either stable or unstable. In the former case the perturbation will
die out over time, and in the latter, the perturbation will blow up
exponentially in time. However, even in the stable case perturbations may get
very large before dying out. The largest amplitude of a perturbed system that
is stable,
\emph{may in fact still grow exponentially in the size of the system}. This
phenomenon is called \emph{flock instability}. Just like ``normal"
instability, flock-instability is an undesirable property, since it makes
large flocks unviable. Flock-instability in arrays of coupled oscillators was
illustrated in~\cite{positionpaper}, and bears similarity to certain phenomena
discovered earlier in fluid mechanics \cite{Trefethen1993,Trefethen1997}. Thus
the first task is to find criteria to identify those systems that are both
stable and flock stable.

We thus need to replace the traditional approach using periodic boundary
conditions by another that we now outline. For those systems that are stable
and flock stable (and \emph{only} for those), we conjecture that for times of
length ${\cal O}(N)$ (where $N$ is the size of the flock) the solutions of the
periodic system behave the same way as in the system with non-trivial
boundary, except near the boundary where additional effects must be taken into
consideration. It turns out that with those constraints the system with
periodic boundary condition behaves like a wave-equation. Since the travel
time of a wave between the leader (agent $0$) and the last agent
(numbered $N$) is proportional to $N$, we can study the dynamics of the
perturbed system for times needed up to a finite number of reflections. Due to
the asymmetry, wave-packages traveling in the positive $\mathbb{R}$ direction
may have a different \emph{signal-velocity} than waves traveling in the
opposite direction. It turns out we can use this effect to achieve either
substantial attenuation or magnification of the traveling wave at the boundary
near agent $N$.

In the present work we extend this analysis from nearest neighbor systems done
in \cite{Cantos2014-1,Cantos2014-2} to next nearest neighbor (NNN) systems,
and in doing that we uncover another new phenomenon. We will see that for
stable and flock stable systems there are still two signal velocities, but
that in contrast with nearest neighbor systems it is possible that they have
the same sign. This means that perturbations can travel (as waves) in only one
direction. As a consequence, they cannot be reflected. This type of transient
has the counter-intuitive characteristic that they travel through the system
in finite time, after which the system finds itself in (almost) perfect
equilibrium.

The paper is organized as follows. In Section~\ref{chap:system} we define the
model of interacting agents. The main line of reasoning of the method is given
in Section~\ref{chap:method}. The stability conditions are given in
Section~\ref{chap:stability}. The description of the stable solutions is
presented in Section~\ref{chap:classification}. This includes the description
of the reflectionless waves on the line, which to the best of our knowledge is
a new result. We include extensive numerical analysis in
Section~\ref{chap:numerics} to back up our theory.

\section{The Equations of Motion of the NNN System}
\label{chap:system}

\begin{figure}[!ht]
\begin{center}
\includegraphics[width=0.75\columnwidth]{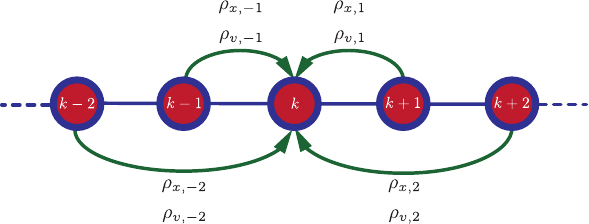}
\end{center}
\caption{{\bf Sketch of information flow.} Available information about
position $\rho_{x,j}$ and velocity $\rho_{v,j}$ weight of nearest $j=k\pm1$
and next nearest $j=k\pm2$ agents for $k$'th agent.}
\label{Fig1}
\end{figure}

We consider a model of an one--dimensional array of linear damped coupled (up
to next nearest neighbor) harmonic oscillators on the line. The oscillators or
agents are numbered from 0 to $N$ from right to left. We impose that the
system is \emph{decentralized}, that is: the agents perceive only information
about other agents that is \emph{relative} to themselves, in this case
relative position and relative velocity. See Figure~\ref{Fig1} for a sketch of
information flow.

The equations of motion of such a system can be written as:
\begin{eqnarray}
\ddot{x}_k=\sum_{j=-2}^2
\left[g_x\rho_{x,j}\left(x_{k+j}-x_k+j\,\Delta \right)
+g_v\rho_{v,j}\left(\dot{x}_{k+j}-\dot{x}_k\right)\right]\,,
\label{eq:eqmo1}
\end{eqnarray}
where $\Delta$ is the desired inter--agent distance and $\rho_{x,j}$
($\rho_{v,j}$) are position (velocity) parameters. The latter are normalized
so that $\rho_{x,0}=\rho_{v,0}=1$. The normalization factors, $g_x$ and
$g_v$, are often called the `gains' in the engineering literature.

The initial conditions we will impose from here on, are as follows. At time
$t\leq0$ the agents are in equilibrium, $x_k=-k\,\Delta$. Then, for $t>0$,
the leader $x_0$ starts moving forward at velocity $v_0$:
\begin{eqnarray*}
\forall t\geq 0\quad x_0(t)=v_0t\,.
\end{eqnarray*}
The leader is not influenced by other agents, although other agents (e.g.,
$k=1$ and $k=2$) are influenced by it. This choice of the leader at the head
of the flock is motivated by applications to traffic situations
(see~\cite{Cantos2014-1,Cantos2014-2}). It is possible to analyze the dynamics
with leaders in different positions or having more than one leader. We will
not pursue this here.

It is convenient to eliminate the constant $\Delta$ from
Equation~(\ref{eq:eqmo1}), using the change of coordinates:
$z_k\equiv x_k+k\,\Delta$. In this notation, the equation of motion of the
flock in $\mathbb{R}$ becomes:
\begin{defn}
The equations of motion of the NNN system with $N>4$ agents, for
$k\in\{1,\cdots N\}$, are:
\begin{eqnarray}
\ddot{z}_{k}=\sum_{j=-2}^{2}
\left(g_{x}\rho_{x,j}z_{k+j}+g_{v}\rho_{v,j}\dot{z}_{k+j}\right)\,.
\label{eq:sys}
\end{eqnarray}
This system is subject to the constraints
\begin{eqnarray}
\rho_{x,0}=\rho_{v,0}=1\,,
\quad\sum_{j=-2}^{2}\rho_{x,j}=\sum_{j=-2}^{2}\rho_{v,j}=0\,,
\label{eq:dec}
\end{eqnarray}
and to the initial conditions:
\begin{eqnarray*}
z_k(0)=0\,,\; \dot z_k(0)=0\,,\mathrm{and}\,\;z_0(t)=v_0 t\,.
\end{eqnarray*}
From now on we denote this system by ${\cal S}_N$. The collection of the
systems $\left\{{\cal S}_N\right\}_{N>4}$ will be denoted by ${\cal S}$.
\label{def:nnn}
\end{defn}

Now we use vector notation and write
$z\equiv(z_{1},\;z_{2},\;z_{3},\;\dots\;z_{N}\,)^{T}$ together with
$\dot{z}\equiv(\dot{z}_{1},\;\dot{z}_{2},
\;\dot{z}_{3},\;\dots\;\dot{z}_{N}\,)^{T}$.
Equation~(\ref{eq:sys}) may be rewritten as a first order system in $2N$ dimensions:
\begin{eqnarray}
\frac{\mathrm{d}}{\mathrm{d}t}\pmatrix{z \cr \dot{z}}=
\pmatrix{0 & I \cr L_{x} & L_{v}}\pmatrix{z \cr \dot{z}}+F(t)
\equiv M_N\pmatrix{z \cr \dot{z}}+F(t)\,,
\label{eq:focs}
\end{eqnarray}
where $L_{x}\,,L_{v}\in\mathbb{R}^{N\times N}$ are matrices - the
\emph{Laplacians} - with standard definition
\begin{eqnarray}
\left(L_{x}z\right)_k=\sum_{j=-2}^{2} \, g_{x}\rho_{x,j}z_{k+j}\,, \quad
\left(L_{v}\dot z\right)_k=\sum_{j=-2}^{2}\,g_{v}\rho_{v,j}\dot z_{k+j}\,,
\label{eq:Laplacians}
\end{eqnarray}
where $F(t)$ is the ``external force" that describes the influence of the
leader with trajectory $z_0(t)=v_0t$ on the acceleration of its immediate
neighbors. It is easy to check that all components are zero except the
$N+1$-st and $N+2$-nd components. The exact form of that external force
depends of the boundary conditions we choose to impose on the system, as we
briefly discuss now.

Note that the equations for $z_1$, $z_{N-1}$, and $z_N$ are subject to
non-trivial boundary conditions (BC), because there are no agents with numbers
$-1$, $N+1$, and $N+2$. So the equations of motion for agents $1$, $N-1$, and
$N$, will have to be modified. Here we will use two sets of BC: {\it fixed
interaction} and {\it fixed mass}. In the case of {\it fixed interaction} BC
the central coefficients, $\rho_{x,0}$ and $\rho_{v,0}$, of the boundary
agents are not equal $1$, instead it is the sum of existing interactions. On
the other hand, in {\it fixed mass} BC we change the interactions and keep the
central $\rho$'s equal to $1$. The details are given in the Appendix~A.

\begin{figure}[!hb]
\begin{center}
\includegraphics[width=1.00\columnwidth]{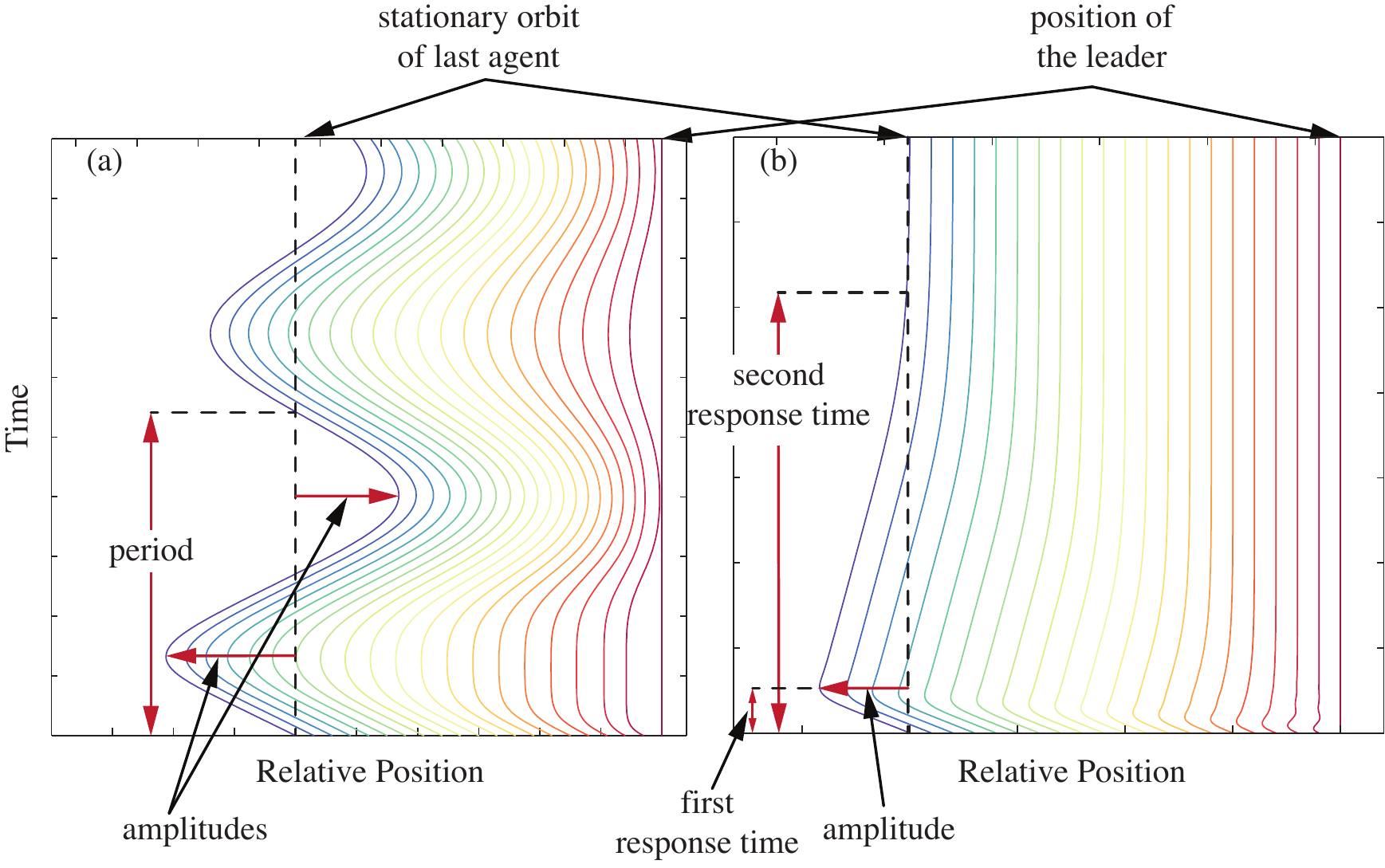}
\end{center}
\caption{{\bf Dynamics of locally coupled arrays.} Sketch of time-dependent
dynamics of locally coupled oscillators on the line (system ${\cal S}_N$) of
(a) Type~I and (b) Type~II (see Section~\ref{chap:classification} for detailed
analysis of these solutions). $x$-axis depicts relative position with respect
to the leader.}
\label{Fig2}
\end{figure}

\emph{Coherent motion} is defined as:
\begin{eqnarray}
y_k(t)=a_0 t +b_0-k\,\Delta\,,
\label{eq:como}
\end{eqnarray}
where $a_0$ and $b_0$ are arbitrary real constants. It is easily checked that
coherent motion is a solution to the differential equations given above. Our
aims are:
\begin{enumerate}
\item To find out for which values of the parameters trajectories the system
is stable: namely, for all $k$, $\lim_{t\rightarrow\infty}|x_k(t)-y_k(t)|=0$
where $y_k$ is given in Equation~(\ref{eq:como}).
\item To find out how fast the stable systems converge to its coherent motion.
\item To determine what is the size of the transient
$\max_{t>0}|x_N(t)-y_N(t)|$ in stable systems.
\end{enumerate}
In the last item we consider only the last (or $N$-th) agent to simplify the
exposition. As an example in Figure~\ref{Fig2} we present a sketch of the
dynamics expected in the {\it stable} system of locally coupled oscillators on
the line. In the figures we plot the positions relative to the leader, i.e.,
$x_k(t)-v_0t$.

\section{Method}
\label{chap:method}

\begin{figure}[!ht]
\begin{center}
\includegraphics[width=1.0\columnwidth]{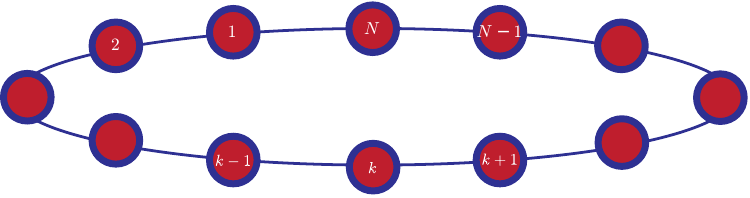}
\end{center}
\caption{{\bf A circular graph.} The numbers indicate the labels of the
agents.}
\label{Fig3}
\end{figure}

The analysis of the system of Definition~\ref{def:nnn} is very difficult
because the Laplacians given in Equation~(\ref{eq:Laplacians}) are not
simultaneously diagonalizable. In order to overcome that we define a system
where the communication structure is not a line graph but a circular graph (
see Figure~\ref{Fig3}). We use it (in Section~\ref{chap:stability}) to
deduce necessary conditions for stability and flock stability on the line, and
(in Section~\ref{chap:classification}) to derive expressions for the signal
velocities.

\begin{defn}
The equations of motion of the system with periodic boundary conditions (PBC)
are:
\begin{eqnarray*}
\ddot{z}_{k}=\sum_{j=-2}^{2}
\left(g_{x}\rho_{x,j}z_{k+j}+g_{v}\rho_{v,j}\dot{z}_{k+j}\right)\,.
\end{eqnarray*}
This system is subject to the constraints
\begin{eqnarray*}
\rho_{x,0}=\rho_{v,0}=1\,,\quad
\sum_{j=-2}^{2}\rho_{x,j}=\sum_{j=-2}^{2}\rho_{v,j}=0\,.
\end{eqnarray*}
Finally, instead of boundary conditions for $z_1$, $z_{N-1}$, and $z_N$, we
set:
\begin{eqnarray*}
\forall j \quad z_{N+j}=z_j\,.
\end{eqnarray*}
From now on we denote this system by ${\cal S}^*_N$. The collection of the
systems $\left\{{\cal S}^*_N\right\}_{N>4}$ will be denoted by ${\cal S}^*$.
\end{defn}
The Laplacians $L^*$ [with the same definition as in
Equation~(\ref{eq:Laplacians})] now become \emph{circulant} matrices and are
therefore diagonalizable by the discrete Fourier transform \cite{Kra}. Let
$w_m$ be the $m$-th eigenvector of $L^*$'s, that is the vector whose $j$-th
component satisfies:
\begin{eqnarray*}
(w_{m})_{j}=e^{\imath\frac{2\pi m}{N}j}\equiv e^{\imath\phi j}\,,
\end{eqnarray*}
with $\phi=2\pi m/N$. We denote the $m$-th eigenvalues of $L^*_x$ by
$\lambda_{x,m}$ and those of $L^*_v$ by $\lambda_{v,m}$. With a slight abuse
of notation we also consider these eigenvalues to be functions
$\lambda_x(\phi)$ and $\lambda_v(\phi)$ of $\phi$ defined above. By using the
$m$-th eigenvector above to calculate $L^*_{x}w_m$ and $L^*_{v}w_m$ it is easy
to show that:
\begin{lemm} The $\lambda$'s are given by
\begin{eqnarray*}
\lambda_{x}(\phi) &=& g_{x}\sum_{j=-2}^{2} \rho_{x,j}\,e^{\imath\phi j}=
g_{x}\sum_{j=0}^{2}
\left[\alpha_{x,j}\,\cos(j\phi)+\imath\beta_{x,j}\,\sin(j\phi)\right]\,,\\
\lambda_{v}(\phi) &=& g_{v}\sum_{j=-2}^{2}\rho_{v,j}\, e^{\imath\phi j}=
g_{v}\sum_{j=0}^{2}
\left[\alpha_{v,j}\,\cos(j\phi)+\imath\beta_{v,j}\,\sin(j\phi)\right]\,.
\end{eqnarray*}
\label{lemm:the lambdas}
\end{lemm}

Here we have used the following convenient notation.
\begin{defn}
Let $\alpha_{x,0}=\alpha_{v,0}=1$ and $\beta_{x,0}=\beta_{v,0}=0$. For $j>0$
we define:
\begin{eqnarray*}
\alpha_{x,j}=\rho_{x,j}+\rho_{x,-j}\,,
\quad\beta_{x,j}=\rho_{x,j}-\rho_{x,-j}\,,\\
\alpha_{v,j}=\rho_{v,j}+\rho_{v,-j}\,,
\quad\beta_{v,j}=\rho_{v,j}-\rho_{v,-j}\,.
\end{eqnarray*}
Note that the sum of the $\alpha$'s equals $0$ by Equation~(\ref{eq:dec}).
\label{defn:new parameters}
\end{defn}

Let us now focus on the eigenvectors and eigenvalues of $M^*_N$ associated
with $w_m$. Denoting the eigenvalues by $\nu_{m,\pm}$, we get:
\begin{eqnarray}
\pmatrix{0 & I \cr L^*_{x} & L^*_{v}}
\pmatrix{w_{m} \cr \nu_{m,\pm}\,w_{m}}=
\nu_{m,\pm}
\pmatrix{w_{m} \cr \nu_{m,\pm}\,w_{m}}\,.
\label{eq:mat}
\end{eqnarray}
Thus the evolution of an arbitrary initial condition is given by:
\begin{eqnarray}
\pmatrix{z(t) \cr \dot z(t)} =
\sum_m\,a_m e^{\nu_{m,-}t}\,\pmatrix{w_{m} \cr \nu_{m,-}\,w_{m}}+
\sum_m\,b_m e^{\nu_{m,+}t}\,\pmatrix{w_{m} \cr \nu_{m,+}\,w_{m}}\,,
\label{eq:evolution}
\end{eqnarray}
where the $a_m$ and $b_m$ are determined by the initial condition at $t=0$.

Next, let us evaluate the second row of Equation~(\ref{eq:mat}) using that
$w_m$ are eigenvectors of $L^*$:
\begin{lemm}
The eigenvalues of ${\cal S}^*_N$ are the roots of the characteristic equation
\begin{eqnarray}
\nu^2-\lambda_v(\phi)\nu-\lambda_x(\phi)=0\,,
\label{eq:char}
\end{eqnarray}
where $\phi=2\pi m/N$. The eigenvalues of ${\cal S}^*$ are a dense subset of
the closed curves $\nu_+:S^1\rightarrow \mathbb{C}$ and
$\nu_-:S^1\rightarrow \mathbb{C}$ defined by Equation~(\ref{eq:char}).
\label{lemm:char}
\end{lemm}

Our treatment follows that of~\cite{Cantos2014-2} where it is conjectured that
(for nearest neighbor systems) a circular system and a system on the line
evolve in a similar manner. The result is that we can analyze the circular
system and apply the conclusions to the systems on the line. We briefly
outline how the evolution of the two systems can be compared.

First we need to remind the reader of the two notions of stability that play a
crucial role in our analysis.
\begin{defn}
For given $N$, the system ${\cal S}_N$ (${\cal S}^*_N$) is asymptotically
stable if, given any initial condition, the trajectories always converge to a
coherent motion and the convergence is exponential in time. For the systems we consider
this is equivalent
to: $M_N$ ($M^*_N$) has one eigenvalue zero with multiplicity $2$, and all other
eigenvalues have real part (strictly) less than $0$. ${\cal S}_N$
(${\cal S}^*_N$) is unstable if at least one eigenvalue has positive real part.
\label{def:asymptstable}
\end{defn}
Flock stability was introduced in~\cite{positionpaper}:
\begin{defn}
The collection ${\cal S}$ is called flock stable if the ${\cal S}_N$ are
asymptotically stable for all $N$ and if
$\mathrm{max}_{t\in\mathbb{R}}|z_{N}(t)|$ grows sub--exponentially in $N$.
\end{defn}
Note that asymptotic stability is different from flock stability. The former
deals with the growth of the response of a single system as $t$ tends to infinity
while $N$ is held fixed, while
the latter deals with the growth of the response of a sequence of systems as
$N$ tends to infinity.

Now we mention the main ideas that allow us to compare the evolution of the
two systems. The first idea is the conjecture that if the system on the circle
is asymptotically unstable, then the system on the line is either
asymptotically unstable or flock unstable. Notice that undamped, symmetric
systems are all marginally stable, and this aspect does not enter the
traditional discussion in the physics context. This gives us necessary
conditions for stability and flock stability of the system on the line.

The second idea involved in this analysis is the principle that, if the system
on the line is stable and flock stable, then the evolution away from the
boundary of the two systems should be the same. As we shall see this means
that for these systems we obtain wave-like behavior with signal velocities
determined by the eigenvalues of the system on the circle
(see~Theorem~\ref{thm:velocities}). This is similar to what is commonly known
in solid state physics as \emph{periodic boundary conditions} (see Chapter 21
in \cite{Ashcroft1976}), though not exactly the same. The difference is that
here we apply principle in more generality than is usual in physics, because
we are considering systems that are not symmetric and not Hamiltonian.

\begin{figure}[!ht]
\begin{center}
\includegraphics[width=0.75\columnwidth]{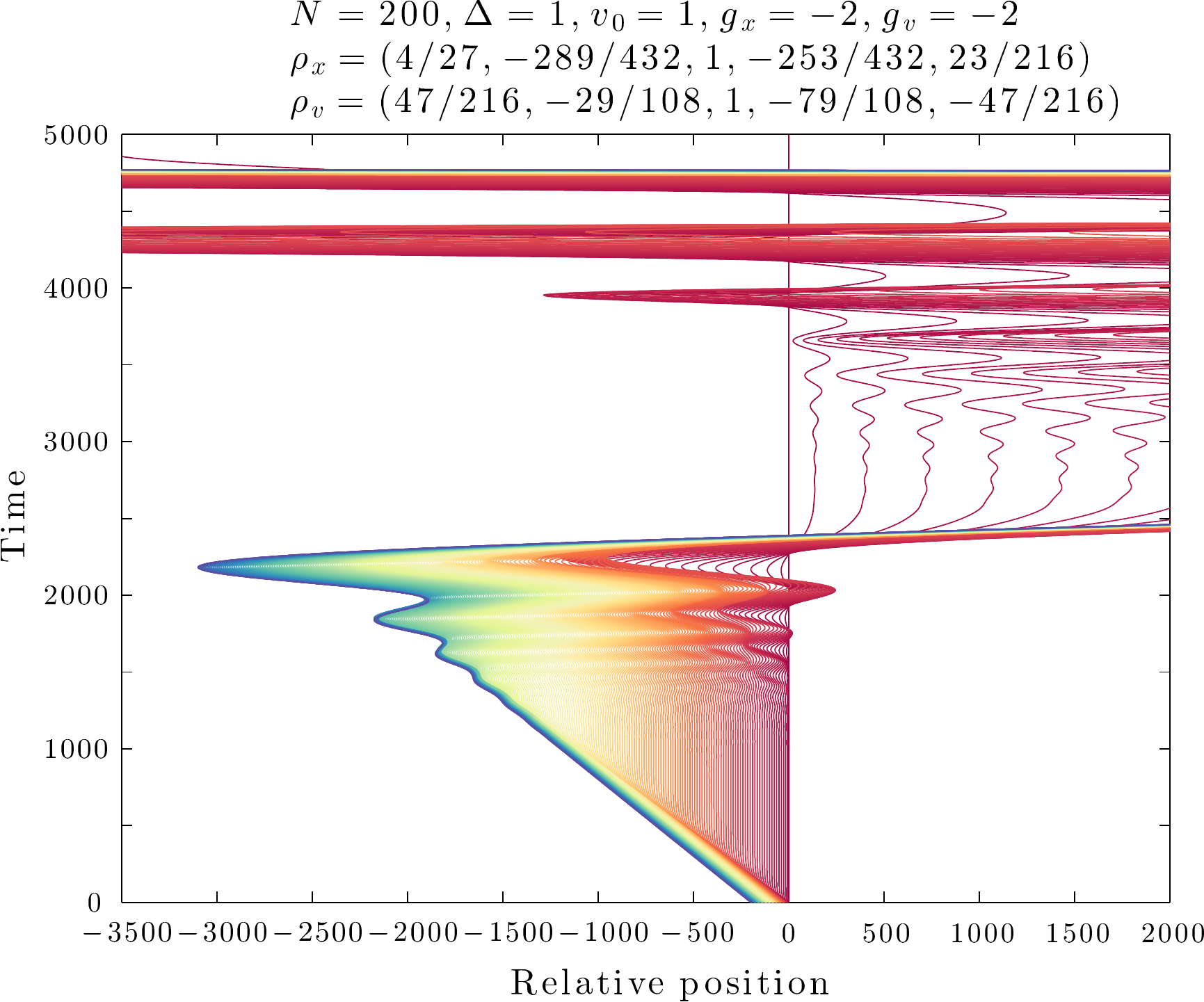}
\end{center}
\caption{{\bf Dynamics of unstable system.} Dynamics of example system
${\cal S}_N$ as calculated for $N=200$, $\Delta=1$, $v_0=1$, $g_x=g_v=-2$,
$\rho_x=(4/27, -289/432, 1, -253/432, 23/216)$,
$\rho_v=(47/216, -29/108, 1, -79/108, -47/216)$ and fixed interaction BC. Each
color represents the orbit of one of the $200$ agents.}
\label{Fig4}
\end{figure}

We know that the reverse of this conjecture is actually \emph{false}:
stability on the circle \emph{does not} imply stability on the line. There are
systems that are stable if periodic boundary conditions are imposed, but have
some eigenvalues with positive real parts when given non-trivial physical
boundary conditions. In Figure~\ref{Fig4} we show a simulation of such a
system on the line. The parameters are given in the Figure. Another example is
given in Section~\ref{chap:classification}. It turns out, perhaps fortunately,
that such \emph{counter examples} are extremely rare.

The third and last idea is that the cumbersome physical boundary conditions of
Appendix~A may be replaced by a single ``free boundary condition" and a single
``fixed boundary condition". This is a great simplification, because the set
of possible all physical boundary conditions form a $16$-parameter set, with
no obvious naturally ``preferred" boundary condition. However because of this
last principle, our conclusions will be independent of the physical boundary
condition. As before, in the traditional physics context, this problem play
little or no role, because presumably the fixed mass BC is the only possible
BC.

In extending the principle of periodic boundary conditions and adding some new
ideas to it, we need to be aware that new phenomena may appear (see
Section~\ref{chap:subseII}) and indeed its validity is not guaranteed nor is
it implied by the validity of the principle in the restricted (symmetric,
undamped) case (nor indeed by the validity in the general nearest neighbor
case). Thus our conclusions need to be checked numerically (see
Section~\ref{chap:numerics}).

\section{Stability Conditions}
\label{chap:stability}

We wish to establish conditions that guarantee that the systems ${\cal S}_N$
on the line is both asymptotically stable and flock stable. Since a direct
verification is too hard or even impossible to perform, we use the conjectures
stated in Section~\ref{chap:method}. According to those, necessary conditions
include the stability of the systems ${\cal S}^*$, a much simpler problem.

Substituting the expressions for the $\lambda$'s in
Lemma~\ref{lemm:the lambdas} into Equation~(\ref{eq:char}), we see that the
eigenvalues of ${\cal S}^*_N$ are the roots of the following equation:
\begin{eqnarray}
\nu^{2}-\nu \, g_{v}\sum_{j=-2}^{2}\rho_{v,j}\,e^{\imath\phi j} -
g_{x}\sum_{j=-2}^{2}\rho_{x,j}\,e^{\imath\phi j}=0
\label{eq:characteristic}
\end{eqnarray}
Note that when $\phi=0$, the characteristic equation becomes $\nu^2=0$. This
gives two zero eigenvalues. These trivial eigenvalues are associated with the
coherent solutions of the system, $z_k=0$ [see also Equation~(\ref{eq:como})].

\begin{lemm}
The following are necessary conditions for ${\cal S}_N^*$ not to have
eigenvalues with positive real part when $N$ is large:
\begin{description}
\item[(i)] $\beta_{x,1}+2\beta_{x,2}=0$,
\item[(ii)] $g_{v}\leq0$,
\item[(iii)] $\alpha_{v,1}\in[-4/3,0]$,
\item[(iv)] $g_{x}\alpha_{x,1}\geq 0$.
\end{description}
\label{lemm:simple conditions}
\end{lemm}
{\bf Proof:} To prove {\bf (i)} notice that the roots of characteristic
Equation~(\ref{eq:char}) are:
\begin{eqnarray}
\nu_{\pm}(\phi)=\frac{1}{2}\left[\lambda_{v}(\phi)
\pm\sqrt{\lambda_{v}(\phi)^2+4\lambda_{x}(\phi)}\right]\,.
\label{eq:roots nu}
\end{eqnarray}
As $\phi=2\pi m/N$ becomes very small, the $\lambda$'s can be approximated by
their first order expansion. From Definition~\ref{defn:new parameters} and
Lemma~\ref{lemm:the lambdas} we obtain:
\begin{eqnarray*}
\lambda_{x}(\phi\to0)\approx\imath g_{x}\phi\sum_{j=0}^{2} j\beta_{x,j}\,,\quad
\lambda_{v}(\phi\to0)\approx\imath g_{v}\phi\sum_{j=0}^{2} j\beta_{v,j}\,.
\end{eqnarray*}
Substituting these into equation for $\nu$, Equation~(\ref{eq:roots nu}), we
see that for small enough $\phi$, the term $\pm\sqrt{4\lambda_{x}(\phi)}$
dominates. Since $\phi$ can be either positive or negative, this has four
branches meeting at the origin at angles of $\pi/2$. Two of these branches
contain eigenvalues with positive real part (for big enough $N$). Therefore,
for $N$ large enough there are $\phi$ such that $\nu_\pm(\phi)$ have negative
real part \emph{unless} $\sum_{j=0}^{2} j\beta_{x,j}=0$.

For condition {\bf (ii)} we note that the mean of the two roots of
Equation~(\ref{eq:roots nu}) is equal $\lambda_v/2$. It follows that we must
require $\Re[\lambda_v(\phi)]\leq0$ for all $\phi\neq 0$. Since the average
$\frac{1}{2\pi}\,\int_{-\pi}^{\pi} \Re[\lambda_v(\phi)]\mathrm{d}\phi$ is
$g_v$, there is a $\phi$ so that $\Re[\lambda_v(\phi)]\geq g_v$. That of
course means that $g_v$ must be non-positive.

For {\bf (iii)} note that $\Re[\lambda_v(\phi)]\leq 0$. Therefore,
$\sum \alpha_{v,j} \cos j\phi\geq 0$. For the NNN system, the constraints on
the $\alpha$'s now give
\begin{eqnarray*}
1+\alpha_{v,1}\cos(\phi)-(1+\alpha_{v,1})\cos(2\phi) \geq 0\,.
\end{eqnarray*}
Since $\cos(2\phi)=2\cos^2\phi-1$, the inequality becomes a quadratic
inequality in $\cos(\phi)$:
\begin{eqnarray*}
-(2+2\alpha_{v,1})\cos^2(\phi)+\alpha_{v,1}\cos \phi+2+\alpha_{v,1}\geq0\,,
\end{eqnarray*}
which factors as:
\begin{eqnarray*}
-\left[(2+2\alpha_{v,1})\cos(\phi)+2+\alpha_{v,1}\right]
\left( \cos(\phi)-1\right)\geq 0\,.
\end{eqnarray*}
By working out three cases, $\alpha_{v,1}<-1$, $\alpha_{v,1}=-1$, and
$\alpha_{v,1}>-1$, the conclusion of {\bf (iii)} may be verified.

Beside $\phi=0$, one other case of Equation~(\ref{eq:char}) is easy, namely
$\phi=\pi$ with the $\lambda$'s as defined in Lemma~\ref{lemm:the lambdas}
\begin{eqnarray*}
\nu^2-\nu g_v\sum_{j= 0}^{2}(-1)^j\alpha_{v,j}-g_x
\sum_{j= 0}^{2}(-1)^j\alpha_{x,j}=0\,.
\end{eqnarray*}
The roots have non-positive real part if and only if both coefficients are
nonnegative. In particular, this implies that last term in the above equation
is $g_x\sum_{j=0}^{2}{(-1)^j\alpha_{x,j}}\leq 0$. From
Definition~\ref{defn:new parameters} we know that
$\sum_{j=0}^{2}{\alpha_{x,j}}=1+\sum_{j=1}^{2}{\alpha_{x,j}}=0$, and as a
consequence $g_x\alpha_{x,1}\geq 0$, which is condition {\bf (iv)}. Similarly,
$g_v \alpha_{v,1}\geq 0$ but this already follows from conditions {\bf (ii)}
and {\bf (iii)}.
\qed

\vskip .1in
Since we are only interested in the parameter values for which the collection
${\cal S}^*$ is not unstable, we use the above
Lemma~\ref{lemm:simple conditions} and Definition~\ref{defn:new parameters} to
eliminate a few parameters from our equations. This is done by eliminating
$\beta_{x,2}$, $\alpha_{x,2}$, and $\alpha_{v,2}$ through the substitution:
\begin{eqnarray*}
\beta_{x,2}=-\frac{1}{2}\beta_{x,1}\,,\quad \alpha_{x,2}=-(1+\alpha_{x,1})\,,
\quad\alpha_{v,2}=-(1+\alpha_{v,1})\,,
\end{eqnarray*}
which we will use from now on.

\begin{prop}
If the collection ${\cal S}^*$ is stable, the low-frequency expansion of
$\nu_\pm(\phi)$ is given by
\begin{eqnarray*}
\nu_{\pm}(\phi) &=& \frac{\imath \phi}{2}\left[
g_{v}(\beta_{v,1}+2\beta_{v,2})\pm\sqrt{g_{v}^{2}
(\beta_{v,1}+2\beta_{v,2})^{2}-2g_{x}(4+3\alpha_{x,1})}\right]\\
&+&\frac{\phi^2}{4}\left[
g_{v}(4+3\alpha_{v,1})\pm\frac{g_{v}^{2}(\beta_{v,1}+2\beta_{v,2})
(4+3\alpha_{v,1})
+2g_{x}\beta_{x,1}}{\sqrt{g_{v}^{2}(\beta_{v,1}+2\beta_{v,2})^{2}-2g_{x}
(4+3\alpha_{x,1})}}\right]\,.
\end{eqnarray*}
\label{prop:expansion of nu}
\end{prop}
{\bf Proof:} One can transcribe the first two terms of the corresponding
expansion given in \cite{Cantos2014-1}, or one can find the result by
substituting power series in $\phi$ in Equation~(\ref{eq:characteristic}) or
Equation~(\ref{eq:roots nu}).
\qed

\vskip .1in
This result immediately implies two other necessary criteria for stability. It
is unclear whether together with the earlier criteria from
Lemma~\ref{lemm:simple conditions} these also constitute a sufficient set of
criteria for the stability of ${\cal S^*}$.

\begin{cory}
The following are necessary conditions for the collection ${\cal S}^*$ to not
be unstable:
\begin{description}
\item[(i)] $g_{v}^{2}(\beta_{v,1}+2\beta_{v,2})^{2}-2g_{x}
(4+3\alpha_{x,1})\geq0$,
\item[(ii)] $g_{v}^{2}g_{x}(4+3\alpha_{v,1})^{2}
(4+3\alpha_{x,1})+2g_{v}^{2}g_{x}(\beta_{v,1}
+2\beta_{v,2})(4+3\alpha_{v,1})\beta_{x,1}+2g_{x}^{2}\beta_{x,1}^{2}\leq0$.
\end{description}
\label{cor:extra conditions}
\end{cory}
{\bf Proof:} If condition {\bf (i)} does not hold, then one branch of the
first order expansion given in Proposition~\ref{prop:expansion of nu} will
have positive real part. Condition {\bf (ii)} corresponds to setting the
argument of $\phi^2$ in Proposition~\ref{prop:expansion of nu} as negative.
\qed

\begin{rema}
We summarize the stability criteria for later use. From Lemma~\ref{lemm:simple conditions} and
Corollary~\ref{cor:extra conditions} we get a list of necessary conditions for system
stability. We added condition vii which was derived in  Corollary~\ref{cor:rhc} 
using Routh--Hurwitz stability criteria (details are given in Appendix~B). 
\begin{description}
\item[{\bf(i)}] $\beta_{x,1}+2\beta_{x,2}=0$,
\item[{\bf(ii)}] $g_{v}\leq0$,
\item[{\bf(iii)}] $\alpha_{v,1}\in[-4/3,0]$,
\item[{\bf(iv)}] $g_{x}\alpha_{x,1}\geq 0$,
\item[{\bf(v)}] $g_{v}^{2}
(\beta_{v,1}+2\beta_{v,2})^{2}-2g_{x}(4+3\alpha_{x,1})\geq0$,
\item[{\bf(vi)}] $g_{v}^{2}g_{x}(4+3\alpha_{v,1})^{2}(4+3\alpha_{x,1})
+2g_{v}^{2}g_{x}(\beta_{v,1}+2\beta_{v,2})(4+3\alpha_{v,1})\beta_{x,1}
+2g_{x}^{2}\beta_{x,1}^{2}\leq0$,
\item[{\bf(vii)}] $g_{x}-g_{v}^2\sum_{j=-2}^{2} \rho_{v,j}^2\le 0$
\end{description}
\label{re:listcond}
\end{rema}

\section{Characterization of Solutions}
\label{chap:classification}

We assume that we start with an initial condition given as
Equation~(\ref{eq:evolution}).
\begin{theo}
Let $K_0>0$ fixed. Suppose the collection ${\cal S^*}$ is stable and that the
initial condition is such that there are $\alpha\in(0,1)$ and $q>0$ such that
$Na_m m^{1+q}$ and $Nb_m m^{1+q}$ are bounded, and $(2-q)\alpha\leq 1$. Then
for large $N$ there are functions $f_+$ and $f_-$ such that the solutions
$z_j(t)$ of ${\cal S}_N^*$ satisfy
\begin{eqnarray*}
\lim_{N\rightarrow\infty}
\sup_{t\in[0,K_0N]}|z_j(t)-v_0t-f_-(j-c_-t)-f_+(j-c_+t)|=0\,.
\end{eqnarray*}
The \emph{signal velocities} $c_{\pm}$ are given by
\begin{eqnarray*}
c_{\pm}=-\frac{1}{2}g_{v}(\beta_{v,1}+2\beta_{v,2})
\pm\frac{1}{2}\sqrt{g_{v}^{2}
(\beta_{v,1}+2\beta_{v,2})^{2}-2g_{x}(4+3\alpha_{x,1})}\,.
\end{eqnarray*}
\label{thm:velocities}
\end{theo}
{\bf Sketch of Proof:}
If ${\cal S^*}$ is stable then Definition~\ref{def:asymptstable} and
Lemma~\ref{lemm:char} imply that the eigenvalues lie on curves bounded away
from the imaginary axes, \emph{except} near $\phi=0$ where we have an
eigenvalue $0$ with multiplicity $2$. The low-frequency expansion of $\nu_\pm$
(Proposition~\ref{prop:expansion of nu} and
Corollary~\ref{cor:extra conditions}) implies that in a neighborhood $I_0$ of
$\phi=0$ we can write
\begin{eqnarray*}
\nu_\pm(\phi)=\imath \phi B_{\pm1}+\phi^2B_{\pm2} +\cdots\,,
\end{eqnarray*}
where $B_{\pm1}$,$B_{\pm2}\in\mathbb{R}$ and furthermore $B_{\pm2}<0$. For $N$
large enough, none of the eigenmodes survive long enough to travel around the
system [$t$ of order ${\cal O}(N)$], \emph{except} those with $2\pi m/N$ in
the neighborhood $I_0$. For these wave-numbers and times scales we may now
neglect dissipation.

We use the initial condition of Equation~(\ref{eq:evolution}) with $b_m=0$.
Neglecting dissipation, the evolution of the $j$-th component can then be
written as
\begin{eqnarray*}
z_j(t)=\sum_m\,a_m \, e^{i\phi B_{-1}t}\, e^{i\phi j}=
\sum_m\,a_m\,e^{i\phi(j+ B_{-1}t)}\,.
\end{eqnarray*}
If we write this as $f_+(j-c_+t)$, we see that $c_+=-B_{-1}$. Similarly by
setting $a_m=0$ (instead of $b_m=0$) one shows that $c_-=-B_{+1}$. The general
case follows by superposition of these two. This yields the asymptotic form of
$z_N(t)$.

To actually prove the remainder indeed tends to zero, one needs the assumption
on the decay of the $a_m$ and $b_m$. This part of the argument is given
in~\cite{Cantos2014-1}.
\qed

\vskip .1in\noindent
\begin{rema}
The signal velocities $c_{-}$ and $c_{+}$ are in units of number of agents per
unit time (not in distance per unit time). A positive velocity means going from
the leader towards the last agent.
\end{rema}

\vskip .1in\noindent
Theorem~\ref{thm:velocities} states that if ${\cal S^*}$ is stable, then for
large $N$ the systems ${\cal S}_N^*$ will evolve like a wave equation. From
the conjectures discussed earlier we conclude that the solutions of
${\cal S}_N$ - for large $N$ - will behave the same way, except near
boundaries. Near the boundaries we apply the appropriate boundary conditions (
see below) to get the final solution. This gives linear growth of the
transients, and that cannot be improved upon.

If these conditions are not met, in particular if ${\cal S}^*$ is unstable,
then the conjectures tell us that ${\cal S}$ is either unstable of flock
unstable. In the first case the coherent motions are unstable solutions, and
in the second, transients are exponential in $N$ before dying out.

\begin{figure}[!ht]
\begin{center}
\includegraphics[width=0.75\columnwidth]{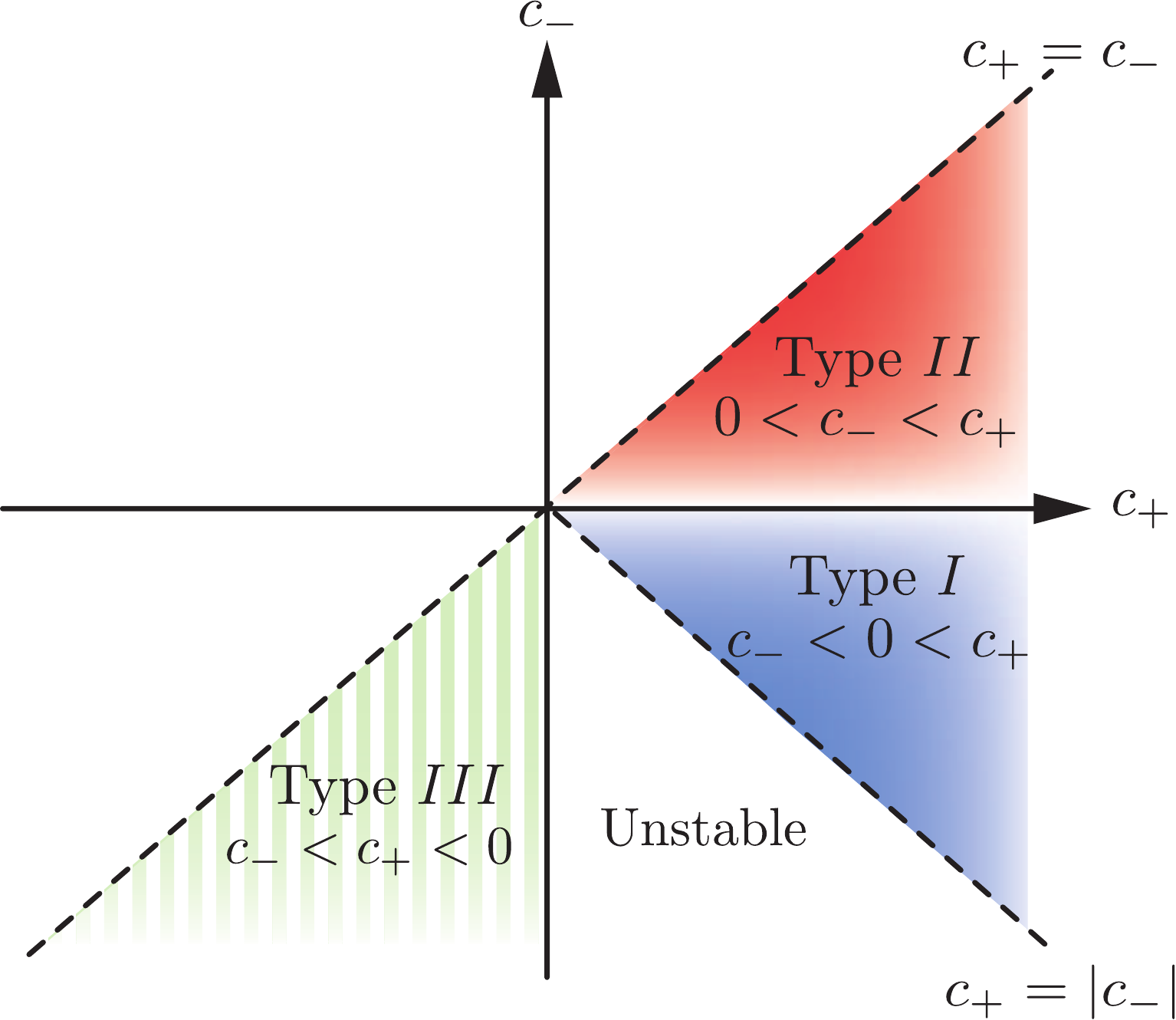}
\end{center}
\caption{{\bf Phase diagram of signal velocities.} Sketch of three type of
solutions. Type~I and Type~II are stable with PBC and also on the line.
Type~III solution is stable only with PBC.}
\label{Fig5}
\end{figure}

It turns out that there are several types of wave-like solutions. These depend
on the signs of the signal velocities $c_{\pm}$ given in
Theorem~\ref{thm:velocities} - see the phase diagram presented in
Figure~\ref{Fig5}. There are, in principle, three types of wave-like solutions. We study these separately.

In our analysis below we ignore cases when $c_{\pm}=0$ or $c_{+}=c_{-}$. These
cases are interesting by themselves, but have properties that make them
undesirable for situations like traffic and other types of flocking. Thus we
do not investigate them here. For example when $c_{-}=0$, distances between
agents don't tend to the desired distance $\Delta$, but rather to some value
that depends on the initial conditions. If $c_{+}=c_{-}$, which only occurs in
Type~II solutions, the velocity of the last agent is unbounded as $N$ tends to
infinity.

\subsection{\bf Type~I: Stable, Flockstable, and $c_{-}<0<c_{+}$}

When $c_{-}<0<c_{+}$ the solutions resemble the traditional damped wave
reflecting between the ends of the flock. The difference in the signal
velocities causes the wave to be damped (or magnified) when it reflects in
agent $N$. These solutions are called Type~I.

For these solutions, it can be shown that the orbit of the last agent [see
Figure~\ref{Fig2}(a)] by the $k$-th \emph{amplitude} $A_k$, the \emph{period}
$T$, and the quotient $|A_{k+1}/A_k|$ which we refer to as the
\emph{attenuation} $\alpha$.

\begin{theo}
Suppose ${\cal S}$ satisfies the conditions of Theorem~\ref{thm:velocities}.
If $c_-<0<c_+$, then for large enough $N$ and at time scales $t={\cal O}(N)$,
the system has Type~I solutions characterized by:
\begin{eqnarray*}
A_{k}=\frac{-v_{0}N}{c_+}\,\left(\frac{c_{-}}{c_{+}}\right)^{k-1}\, ,\quad
\alpha=\left|\frac{c_{-}}{c_{+}}\right|\, , \quad
T=2N\left(\frac{1}{c_{+}}-\frac{1}{c_{-}}\right)\,,
\end{eqnarray*}
where $A_{k}$, $\alpha$, and $T$ are defined above, and $c_\pm$ as in
Theorem~\ref{thm:velocities}.
\label{theo:tyI}
\end{theo}

The proof is essentially that of \cite{Cantos2014-1,Cantos2014-2} and relies
on two insights. The first is that the high frequencies die out fast, so that
we only need to consider low frequencies (as in the proof of
Theorem~\ref{thm:velocities}). The second is that we replace the physical
boundary conditions in ${\cal S}_N$ by new boundary conditions to get the
correct reflections at the ends, namely a \emph{fixed} boundary condition at
the leader's end and a \emph{free} boundary condition at the other end:
\begin{eqnarray*}
z_0(t)=0 \quad \mathrm{and} \quad z_N(t)-z_{N-1}(t)=0\,.
\end{eqnarray*}
Because for large $N$ only low frequencies survive, these conditions can be
replaced by
\begin{eqnarray}
z_0(t)=0 \quad \mathrm{and} \quad
\frac{\partial}{\partial k}\,z_k(t)\big|_{k=N}=0\,.
\label{eq:oldBC}
\end{eqnarray}
That means that near the leader, a pulse reflects (with opposite sign), and
near the free boundary, the traveling pulse is reflected with the same sign
and with amplitude multiplied by a factor $|c_-/c_+|$. The details are written
out in \cite{Cantos2014-2}.

In order to get strong damping to minimize transients, we want
$|c_{-}|<c_{+}$. This means that in the velocity Laplacian, more emphasis
should be placed on the upstream (lower labels) information. Such system have
\emph{asymmetric} interactions.
\begin{cory}
Suppose ${\cal S}$ is asymptotically stable and flock stable. ${\cal S}_N$ has
solutions of Type~I with $|c_-|<c_+$ if: \\
\emph{(i)} $g_v\left(\beta_{v,1}+2\beta_{v,2}\right)<0$ and \\
\emph{(ii)}$g_x\left(4+3\alpha_{x,1}\right)<0$.
\label{cor:attenuation}
\end{cory}
{\bf Proof:} If ${\cal S}$ is asymptotically stable and flock stable, then all
${\cal S}^*_N$ are stable (by our conjectures). The conditions $c_{-}<0<c_{+}$
and $|c_-|<c_+$ imply that $c_-+c_+>0$. This implies \emph{(i)}. Statement
\emph{(i)} together with $c_-<0$ implies \emph{(ii)}.
\qed

\begin{figure}[!ht]
\begin{center}
\includegraphics[width=0.75\columnwidth]{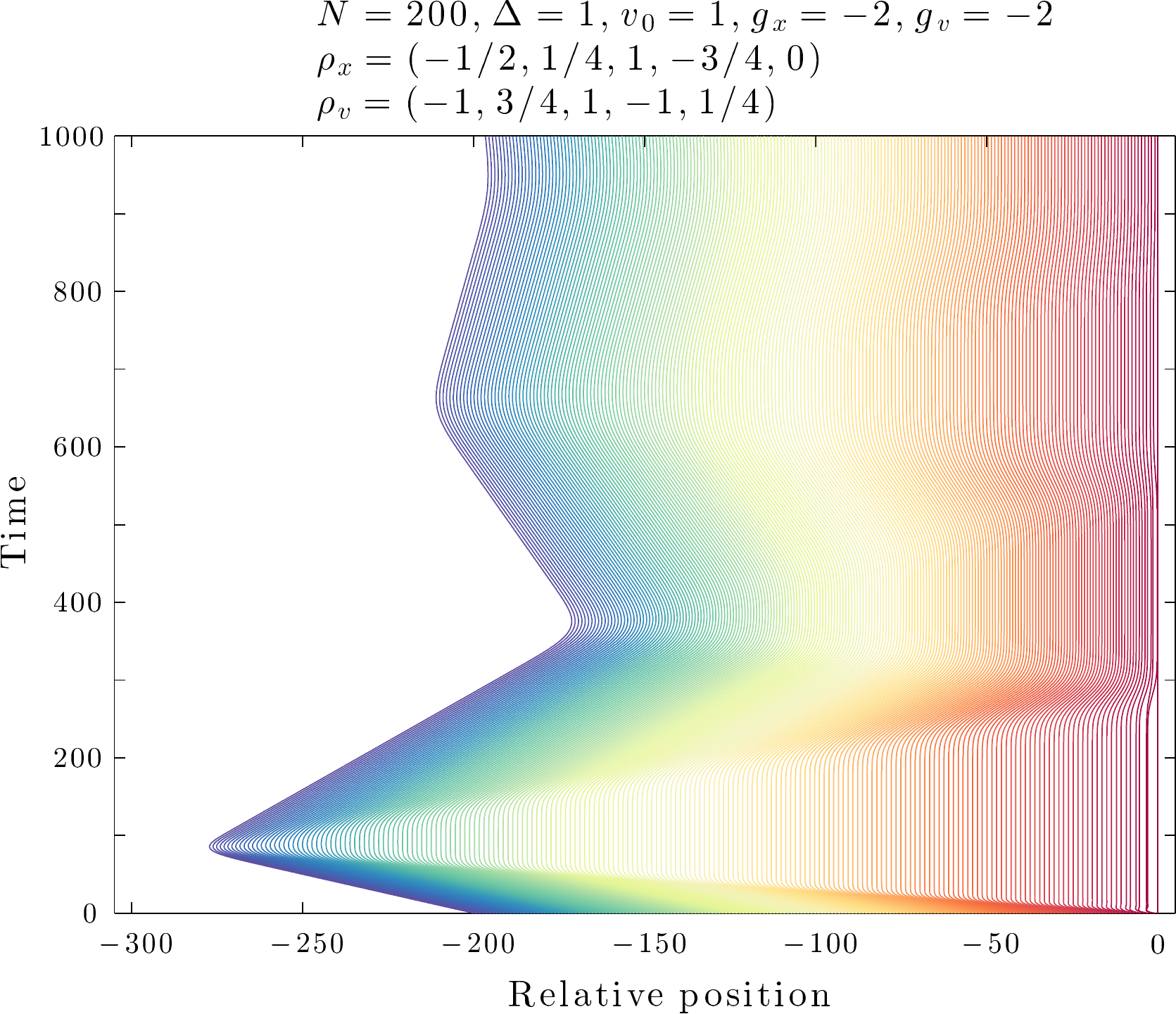}
\end{center}
\caption{{\bf Dynamics of Type~I solution.} Dynamics of example system
${\cal S}_N$ as calculated for $N=200$, $\Delta=1$, $v_0=1$, $g_x=g_v=-2$,
$\rho_x=(-0.5,0.25,1,-0.75,0)$, $\rho_v=(-1,0.75,1,-1,0.25)$ and fixed
interaction BC. Each color represent the orbit of one of the $200$ agents.}
\label{Fig6}
\end{figure}

\vskip .1in\noindent
In Figure~\ref{Fig6} (parameters are given in the figure) we present typical
dynamics of Type~I stable system ${\cal S}_N$. The characteristics predicted
from Theorem~\ref{theo:tyI} are $A_{1}=80\,,\alpha=0.4\,,T=560$, and from the
simulation we measured $A_{1}=77.2$, $\alpha=0.377$, $T=568$.

\subsection{\bf Type~II: Stable, Flockstable, and $0<c_{-}<c_{+}$}
\label{chap:subseII}

When $0<c_{-}<c_{+}$, that is the signal velocities are both positive, the
wave cannot be reflected, because it cannot move with negative velocity. We
denote these solutions as \emph{Type II} or \emph{reflectionless waves}. It
was proved \cite{Cantos2014-2} that such solutions cannot occur with
nearest neighbor interactions.

Since both signal velocities are positive, there is no reflection possible at
$k=N$ agent. Thus the boundary condition at $k=N$ is useless, and we need
another boundary condition. We replace Equation~(\ref{eq:oldBC}) by the
somewhat counter-intuitive condition:
\begin{eqnarray*}
z_0(t)=0 \quad \mathrm{and} \quad z_1(t)-z_0(t)=0\,.
\end{eqnarray*}
As before for large $N$ only low frequencies survive, and so these conditions
can be replaced by
\begin{eqnarray}
z_0(t)=0 \quad \mathrm{and} \quad
\frac{\partial}{\partial k}\,z_k(t)\big|_{k=0}=0\,.
\label{eq:newBC}
\end{eqnarray}
Thus we have both a free and a fixed boundary condition at the leader's end.
For Type~II, the orbit of the last agent [see Figure~\ref{Fig2}(b)] can be
characterized by the \emph{amplitude} $A$, the \emph{first response time}
$T_1$ and the \emph{second response time} $T_2$.

\begin{theo}
Suppose ${\cal S}$ satisfies the conditions of Theorem~\ref{thm:velocities}.
If $0<c_-<c_+$, then for large enough $N$ and at time scales $t={\cal O}(N)$,
the system has Type~II solutions characterized by:
\begin{eqnarray}
A=\frac{-v_{0}N}{c_{+}}\,, \quad
T_{1}=\frac{N}{c_{+}}\,,\quad T_{2}=\frac{N}{c_{-}}\,,
\end{eqnarray}
where $A$, $T_{1}$, and $T_{2}$ are as above, and $c_\pm$ as in
Theorem~\ref{thm:velocities}.
\label{theo:tyII}
\end{theo}
{\bf Proof:} $T_1$ and $T_2$ are the (positive) times at which $z_N(t)-z_0(t)$
changes velocity. These can be deduced from a Proposition whose reasoning is
different enough from earlier work, that we include a sketch of the proof in
Appendix~C. $A=T_1v_0$ is the distance traveled by the leader in the time
interval $[0,T_1)$.
\qed

\begin{cory}
Suppose ${\cal S}$ is asymptotically stable and flock stable. ${\cal S}_N$ has
solutions of Type~II (both velocities positive) if: \\
\emph{(i)} $g_v\left(\beta_{v,1}+2\beta_{v,2}\right)<0$ and \\
\emph{(ii)} $0<2g_x\left(4+3\alpha_{x,1}\right)<g_v^2
\left(\beta_{v,1}+2\beta_{v,2}\right)^2$.
\end{cory}
{\bf Proof:} Similar to the proof of Corollary~\ref{cor:attenuation}.
\qed

\begin{figure}[!ht]
\begin{center}
\includegraphics[width=0.75\columnwidth]{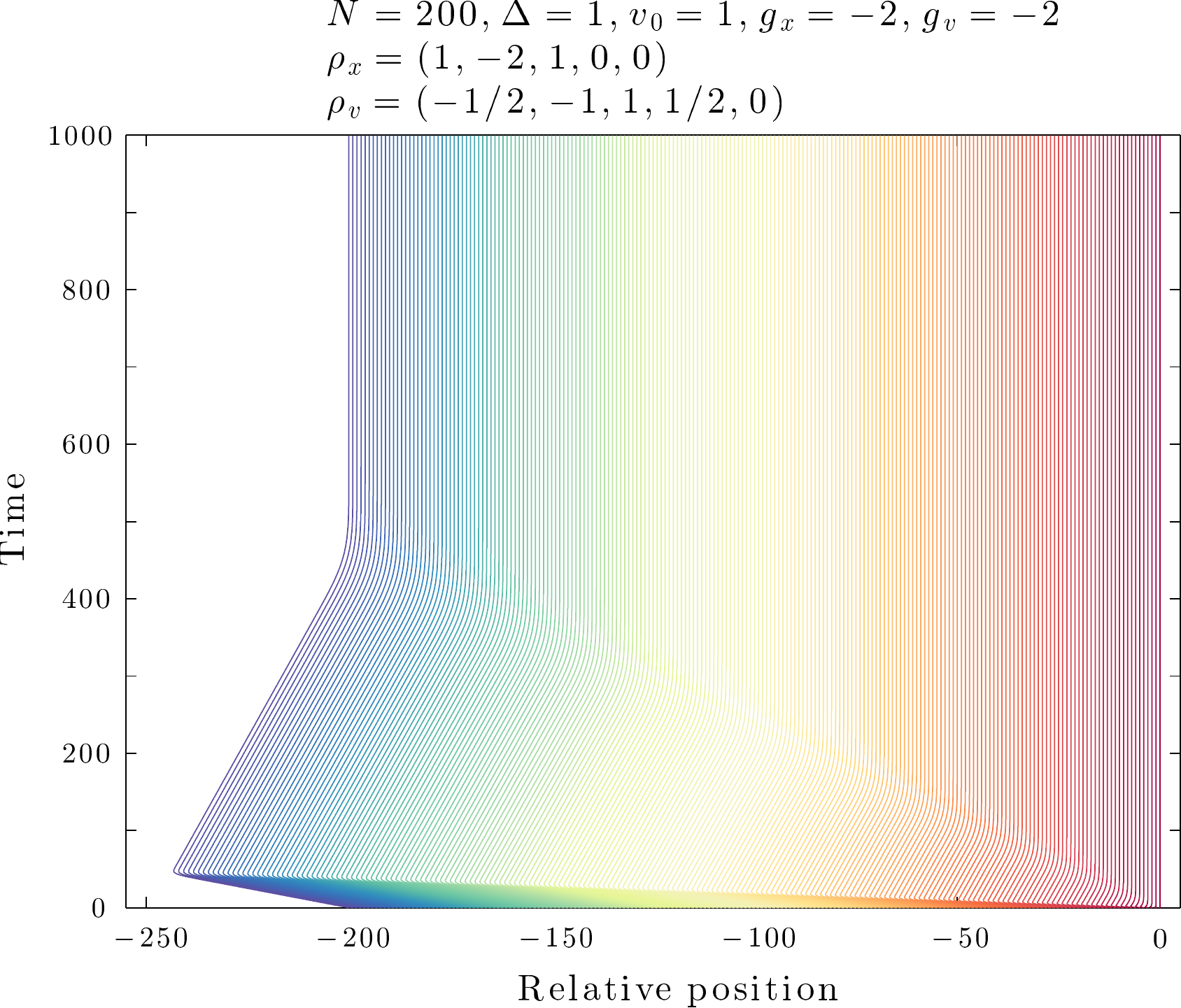}
\end{center}
\caption{{\bf Dynamics of Type~II solution.} Dynamics of example system
${\cal S}_N$ as calculated for $N=200$, $\Delta=1$, $v_0=1$, $g_x=g_v=-2$,
$\rho_x=(1,-2,1,0,0)$, $\rho_v=(-0.5,-1,1,0.5,0)$ and fixed interaction BC.
Each color represent the orbit of one of the $200$ agents.}
\label{Fig7}
\end{figure}

\vskip .1in\noindent
In Figure~\ref{Fig7} we present typical dynamics of Type~II stable system (
parameters given in the figure). The characteristics predicted from
Theorem~\ref{theo:tyII} are $A=43.845$, $T_{1}=43.845$, $T_{2}=456.16$, and
from the simulation we measured: $A=43.182$, $T_{1}=43.182$, $T_{2}=453.95$.
From the figure it seems that a start signal traveling with velocity $c_+$ and
a stop signal traveling with velocity $c_-$ travel from the leader towards the
last agent. A striking aspect of this type of solution is that very briefly
after the second response time, the trajectory of the last agent is (almost)
exactly in its equilibrium position. Dynamics within such a system can be
described as a traveling wave-package which does not reflect in the boundary
of the system.

\subsection{\bf Type~III: $c_{-}<c_{+}<0$}
\label{chap:subseIII}

Finally, when $c_{-}<c_{+}<0$, the perturbation which in our set-up starts at
the leader, cannot be transmitted to the flock, because only negative signal
velocities are available. Thus the system ``finds" another non wave-like
solution which has very large amplitudes. The only reason for listing this
solution in this work at all, is that the system is stable and \emph{does}
have wave-like solutions with negative signal velocities. We call these
solutions Type~III. As with Type~II, these solutions cannot occur with only
nearest neighbor interactions.

\begin{cory}
Suppose ${\cal S}^*$ is asymptotically stable. ${\cal S}_N$ has solutions of
Type~III (both velocities negative) if: \\
\emph{(i)} $g_v\left(\beta_{v,1}+2\beta_{v,2}\right)>0$ and \\
\emph{(ii)} $2g_x\left(4+3\alpha_{x,1}\right)>0$.
\end{cory}
{\bf Proof:} Similar to the proof of Corollary~\ref{cor:attenuation}.
\qed

\begin{figure}[!ht]
\begin{center}
\includegraphics[width=1.0\columnwidth]{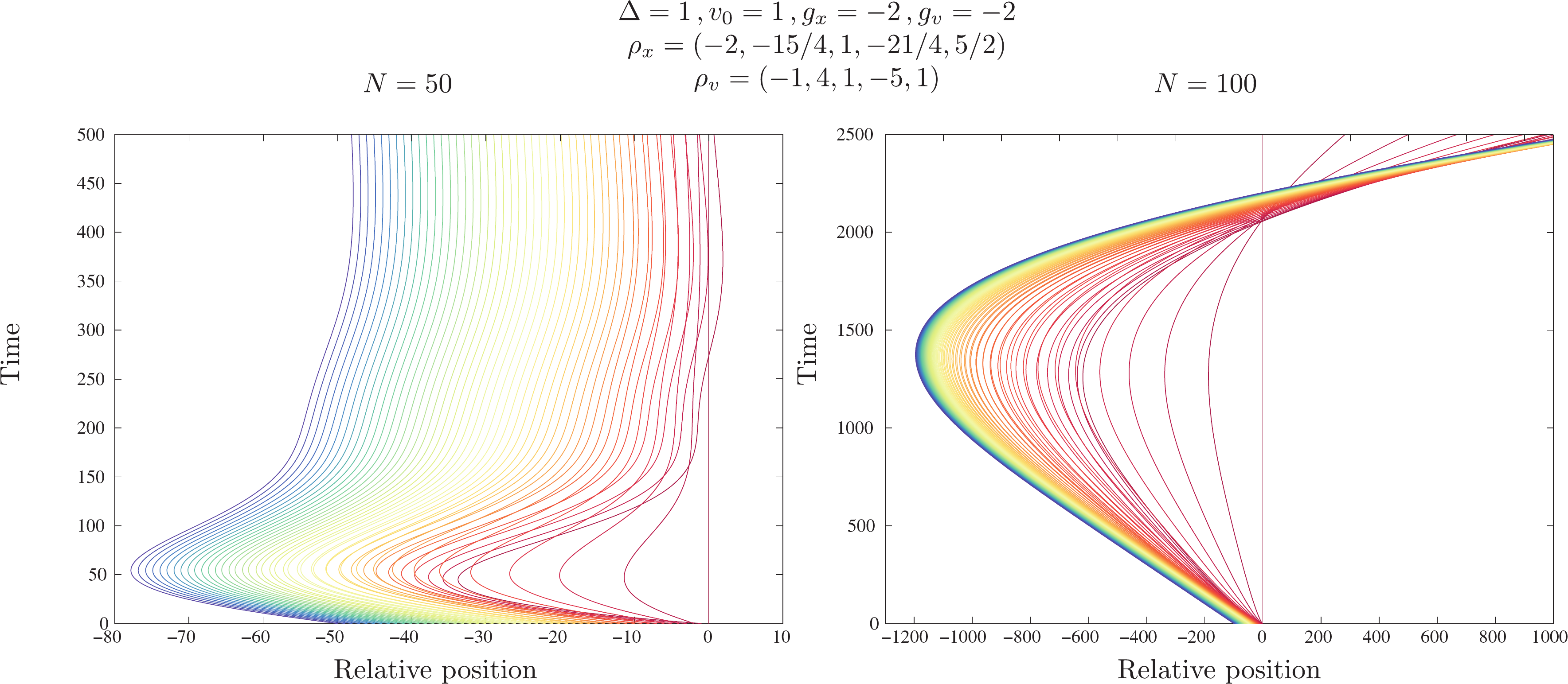}
\end{center}
\caption{{\bf Dynamics of Type~III solution.} Dynamics of example system
${\cal S}_N$ as calculated for $\Delta=1$, $v_0=1$, $g_x=g_v=-2$,
$\rho_x=(-2,-15/4,1,-21/4,5/2)$, $\rho_v=(-1,4,1,-5,1)$, fixed interaction BC
$N=50$ (left panel) and $N=100$ (right panel). Each color represent the orbit
of one of the $50$ and $100$ agents, respectively.}
\label{Fig8}
\end{figure}

\vskip .1in
Within such a setup, on short time scales, the leader simply starts and other
agents do not follow him. On time-scales larger than ${\cal O}(N)$, other
phenomena may take place. Thus amplitudes will grow faster than ${\cal O}(N)$,
and the system is flock unstable or even asymptotically unstable. However, due
to the complicated nature of the stability conditions, we do not have a proof
of this. In Figure~\ref{Fig8} we present a simulation of such a system (
parameters given in the figure). Notice that the amplitudes do not grow
linearly with system size.

\section{Numerical Tests}
\label{chap:numerics}

As we saw in Section~\ref{chap:classification}, measured values of certain
characteristics presented for $N=200$ differ slightly from the predicted ones,
given by Theorem~\ref{theo:tyI} and Theorem~\ref{theo:tyII}. This is expected,
since our predictions are valid for $N\to\infty$. In order to test our
conclusions we performed extensive numerical calculations. We outline our procedure.

First we fixed
$g_x=g_v=-2$ and defined a set $P$ of about $8.6*10^7$ ten-tuples\\
$(g_x,g_v,\rho_{x,-2},\rho_{x,-1},\rho_{x,1},\rho_{x,2},\rho_{v,-2},\rho_{v,-1},\rho_{v,1},\rho_{v,1},\rho_{v,2})$
satisfying Equation~(\ref{eq:dec}). We call these ten-tuples
\emph{configurations}. From this set of configurations we then selected the set
$P_C$ that satisfy all the criteria in Remark~\ref{re:listcond}. For Type~I solutions
we impose an additional constraint, namely: $|c_{-}|<c_{+}$ (see Corollary 
\ref{cor:attenuation}). Next, from
the same ten-tuples of configurations we created the set $P_{S,N}$ that satisfy
Definition~\ref{def:asymptstable} for given $N$. It turns out that for $N$ large
enough these sets were identical: $P_C=P_{S,N}$ (in our case we had to go up to
$N=800$ for a few systems). This strongly suggests that indeed the criteria in
Remark~\ref{re:listcond} (plus Corollary \ref{cor:attenuation}) are a very good indicator of asymptotic stability of
the system on the circle. 

In order to decrease computation time
for large $N$, we imposed a further constraint on $P_S$ that selected $500$
configurations of Type~I and $500$ configurations of Type~II. The constraints
were for the period, namely $T\lesssim {\cal O}(10N)$ (Type~I), and for the
second response time, namely $T_2\lesssim {\cal O}(10N)$ (Type~II). We ran each
of these configurations for $N\in\{25\cdot2^n\}_{n=0}^{n=11}$, that is: for $N$
varying from $25$ to roughly $52,000$. We measure the characteristics directly
from numerical simulations and compare them with predictions of
Theorem~\ref{theo:tyI} and Theorem~\ref{theo:tyII}. For the numerical work we
used the ordinary differential equation solver of the Boost library
\cite{boostlink,boostlib} in a parallel computing environment.

\begin{figure}[!ht]
\begin{center}
\includegraphics[width=1.0\columnwidth]{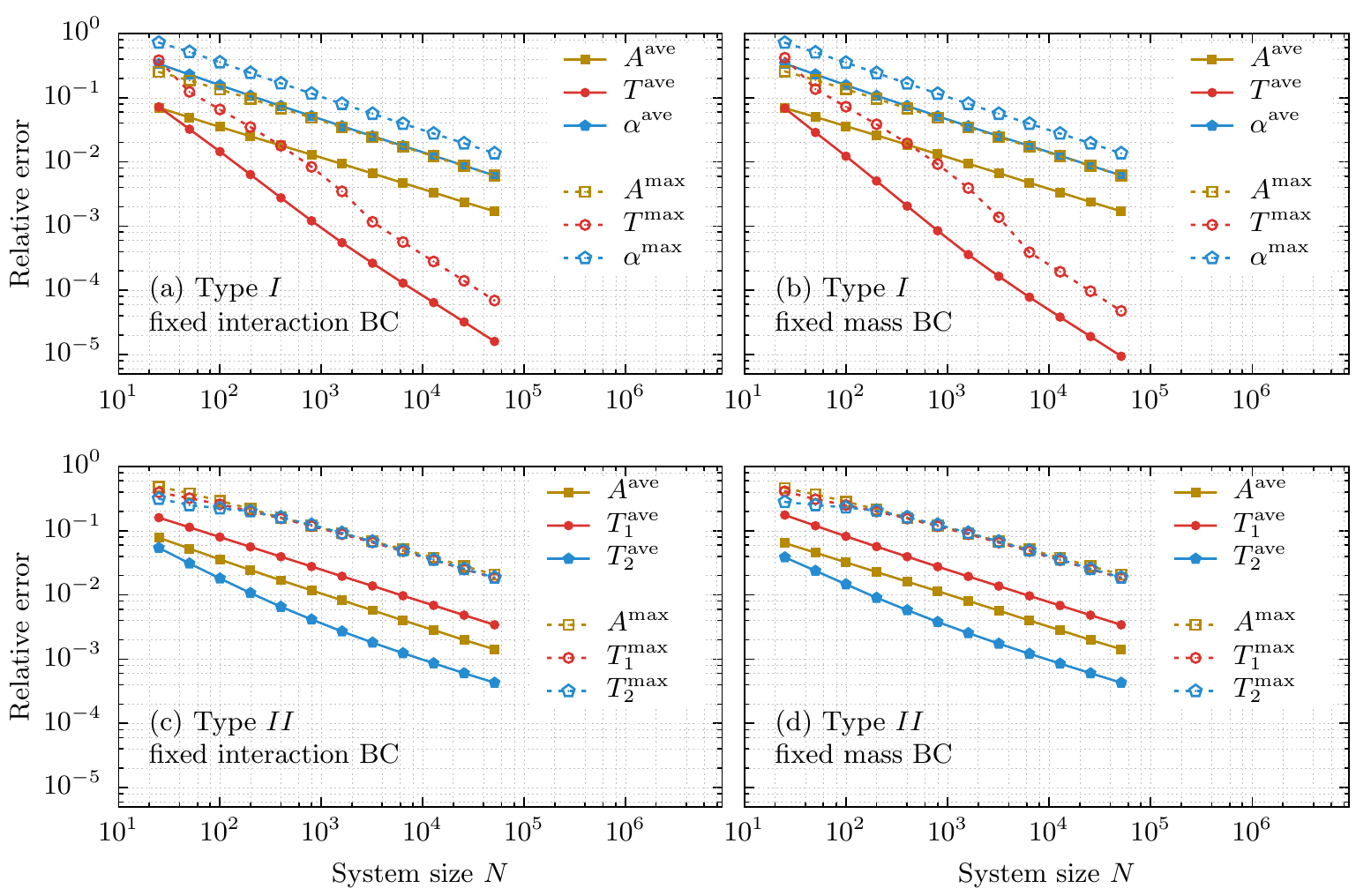}
\end{center}
\caption{{\bf Relative error size scaling.} Size $N$ dependence of average and
maximal relative error of Type~I and Type~II solutions for two different
boundary condition as calculated for $N=25,\dots,51200$ agents. Notice that
the plot has $\log$--$\log$ scale, therefore slope corresponds to the power of
the decay.}
\label{Fig9}
\end{figure}

In Figure~\ref{Fig9} we present the
\emph{relative error}=$|measured-predicted|/|predicted|$ of the following
quantities: for Type~I solutions, the first amplitude $A_1$, the period $T$,
and the attenuation $\alpha$, and for Type~II solutions, the amplitude $A$ and
the first and second response times $T_1$ and $T_2$. We plot both the error
average (for $500$ measurements/configurations) and the worst (largest) error.
We repeated this experiment for two different types of physical boundary
conditions (denoted by {\it fixed interaction} and {\it fixed mass}
(Appendix~A) to make sure that these did not make a difference.

As is clearly visible in Figure~\ref{Fig9}, the relative errors decrease as
$N$ grows, as is predicted by the theory. Our numerical analysis is consistent
with the statement that - with the exception of period $T$ for Type~I orbit -
the error decreases as ${\cal O}(1/\sqrt{N})$. The error in the period $T$
(for Type~I) appears to decrease as ${\cal O}(1/N)$.

\begin{acknowledgements}
We acknowledge support by the European Union's Seventh Framework Program
FP7-REGPOT-2012-2013-1 under grant agreement no.~316165.
\end{acknowledgements}

\section*{Appendix A: The Physical Boundary Conditions}

We will introduce two sets of boundary conditions for ${\cal S}_N$ (the system
on the line). We performed numerics with both types of boundary conditions in
order to support our conclusion that for stable and flock stable systems the
trajectories are independent of the physical boundary conditions.

Let ${\cal S}_N$ be the system in Definition~\ref{def:nnn}. In decentralized
systems the row sum of the Laplacians equals $0$, that is:
$\sum_{j}\rho_{x,j}=\sum_{j}\rho_{v,j}=0$. This implies that for the system
${\cal S}_N$, the equations of agents $k=1$, $N-1$, and $N$ have to be
modified. In the case of fixed interaction BC the masses, $\rho_{x,0}$ and
$\rho_{v,0}$, of the agent are not equal $1$, instead it is the sum of
existing interactions. On the other hand, in fixed mass BC we change the
interactions of existing agents and keep the central $\rho_{x,0}$ and
$\rho_{v,0}$ equal to $1$. Here are the details:
\begin{defn}
\begin{description}
\item[(i)] fixed interaction BC:
\begin{eqnarray*}
\ddot{z}_{1}&=&\left(g_{x}\rho_{x,-1}z_{0}+g_{v}\rho_{v,-1}\dot{z}_{0}\right)\\
&-&\left[g_{x}\left(\rho_{x,-1}+\rho_{x,1}+\rho_{x,2}\right)z_{1}
+g_{v}\left(\rho_{v,-1}+\rho_{v,1}+\rho_{v,2}\right)\dot{z}_{1}\right]\\
&+&\sum_{j=1}^{2}
\left(g_{x}\rho_{x,j}z_{1+j}+g_{v}\rho_{v,j}\dot{z}_{1+j}\right)\,,\\
\ddot{z}_{N-1}&=&\sum_{j=-2}^{-1}\left(g_{x}\rho_{x,j}z_{N-1+j}
+g_{v}\rho_{v,j}\dot{z}_{N-1+j}\right)\\
&-&\left[g_{x}\left(\rho_{x,-2}+\rho_{x,-1}+\rho_{x,1}\right)z_{N-1}
+g_{v}\left(\rho_{v,-2}+\rho_{v,-1}+\rho_{v,1}\right)\dot{z}_{N-1}\right]\\
&+&\left(g_{x}\rho_{x,1}z_{N}+g_{v}\rho_{v,1}\dot{z}_{N}\right)\,,\\
\ddot{z}_{N}&=&\sum_{j=-2}^{-1}\left(g_{x}\rho_{x,j}z_{N+j}
+g_{v}\rho_{v,j}\dot{z}_{N+j}\right)\\
&-&\left[g_{x}\left(\rho_{x,-2}+\rho_{x,-1}\right)z_{N}
+g_{v}\left(\rho_{v,-2}+\rho_{v,-1}\right)\dot{z}_{N}\right]\,.
\end{eqnarray*}
\item[(ii)] fixed mass BC:
\begin{eqnarray*}
\ddot{z}_{1}&=&\left[g_{x}\left(\rho_{x,-2}+\rho_{x,-1}\right)z_{0}
+g_{v}\left(\rho_{v,-2}+\rho_{v,-1}\right)\dot{z}_{0}\right]\\
&+&\sum_{j=0}^{2}
\left(g_{x}\rho_{x,j}z_{1+j}+g_{v}\rho_{v,j}\dot{z}_{1+j}\right) \,,\\
\ddot{z}_{N-1}&=&\sum_{j=-2}^{0}\left(g_{x}\rho_{x,j}z_{N-1+j}
+g_{v}\rho_{v,j}\dot{z}_{N-1+j}\right)\\
&+&\left[g_{x}\left(\rho_{x,1}+\rho_{x,2}\right)z_{N}
+g_{v}\left(\rho_{v,1}+\rho_{v,2}\right)\dot{z}_{N}\right]\,,\\
\ddot{z}_{N}&=&\sum_{j=-2}^{0}\left(g_{x}\rho_{x,j}z_{N+j}
+g_{v}\rho_{v,j}\dot{z}_{N+j}\right)\\
&+&\left[g_{x}\left(\rho_{x,1}+\rho_{x,2}\right)z_{N}
+g_{v}\left(\rho_{v,1}+\rho_{x,2}\right)\dot{z}_{N}\right]\,.
\end{eqnarray*}
\end{description}
\end{defn}

If we use vector notation, the influence of leader on agents $1$ and $2$ is
formulated as an external force. Write
$z\equiv(z_{1},\;z_{2},\;z_{3},\;\dots\;z_{N}\,)^{T}$ and
$\dot{z}\equiv(\dot{z}_{1},\;\dot{z}_{2},\;\dot{z}_{3},\;\dots\;\dot{z}_{N}\,)^{T}$.
The equation of motion can be rewritten as a first order system in $2N$
dimensions:
\begin{eqnarray*}
\frac{\mathrm{d}}{\mathrm{d}t}\pmatrix{z \cr \dot{z}}=
\pmatrix{0 & I \cr L_{x} & L_{v}}\pmatrix{z \cr \dot{z}}+F(t)
\equiv M_N\pmatrix{z \cr \dot{z}}+F(t)\,,
\end{eqnarray*}
Those terms in the full equation of motion that contain $z_0$ or $\dot z_0$
are written as external force. So all components of the external force $F$ are
zero except the $N+1$-st and $N+2$-nd. These two components are given by:
\begin{eqnarray*}
\pmatrix{g_{x}\rho_{x,-1}z_{0}+g_{v}\rho_{v,-1}\dot{z}_{0} \cr
g_{x}\rho_{x,-2}z_{0}+g_{v}\rho_{v,-2}\dot{z}_{0}}\,,
\end{eqnarray*}
if we impose fixed interactions BC, and
\begin{eqnarray*}
\pmatrix{g_{x}\left(\rho_{x,-2}+\rho_{x,-1}\right)z_{0}+
g_{v}\left(\rho_{v,-2}+\rho_{v,-1}\right)\dot{z}_{0} \cr
g_{x}\rho_{x,-2}z_{0}+g_{v}\rho_{v,-2}\dot{z}_{0}}\,,
\end{eqnarray*}
if we impose fixed mass BC.

\section*{Appendix B: The Routh-Hurwitz Stability Criteria}

The Routh-Hurwitz criterion is a standard strategy to derive a concise set of
conditions that is equivalent to the condition that all roots of a given
polynomial have negative real parts. In various systems similar to the ones
discussed here, this criterion gives good results
\cite{Cantos2014-1,Herman2015}. In our current case the resulting equations
are too complicated to give us much information and we only get one more
necessary condition for stability that we can use, namely
Corollary~\ref{cor:rhc}. Our discussion is based on Chapter 15, Sections 6, 8,
and 13 of Ref.~\cite{Gantmacher2000}, where more details can be found.

\begin{theo}
(Routh-Hurwitz) Assume that the determinants given below are nonzero. Given a
real polynomial $R=x^4+a_3x^3+a_2x^2+a_1x+a_0$, all roots of $R$ have negative
real part if and only if all determinants of the upper-left submatrices (the
leading principal minors) of:
\begin{eqnarray*}
A_{4}\equiv
\pmatrix{a_{3} & a_{1} & 0 & 0 \cr 1 & a_{2} & a_{0} & 0 \cr
0 & a_{3} & a_{1} & 0 \cr 0 & 1 & a_{2} & a_{0}}\,,
\end{eqnarray*}
are positive. That is: $a_3>0$, $a_0>0$, $a_3a_2-a_1>0$, and
$a_3a_2a_1-a_3^2a_0-a_1^2>0$.
\label{thm:Hurwitz}
\end{theo}
An equivalent but less well--known set of conditions is given in the following:
\begin{theo}
(Li\'{e}nard-Chipart) Assume that the determinants in
Theorem~\ref{thm:Hurwitz} are nonzero. Given a real polynomial
$R=x^4+a_3x^3+a_2x^2+a_1x+a_0$, all roots of $R$ have negative real part if
and only if $a_3>0$, $a_2>0$, $a_0>0$, and $a_3a_2a_1-a_3^2a_0-a_1^2>0$.
\end{theo}

The characteristic polynomial $Q$ of Equation~(\ref{eq:char}) can be turned into
a polynomial with real coefficients
\begin{eqnarray*}
R=QQ^{*} &\equiv& \nu^4-2\Re(\lambda_v)\nu^3+
\left[|\lambda_v|^2-2\Re(\lambda_x)\right]\nu^2\\
&+& 2\left[\Re(\lambda_x)\Re(\lambda_v)+
\Im(\lambda_x)\Im(\lambda_v)\right]\nu +|\lambda_x|^2\,,
\end{eqnarray*}
by taking its product with its complex conjugate. Clearly, all roots of $Q$
have negative real part if and only if the same is true for $R$. Notice that
in each of the two criteria, one of the equations is trivially satisfied,
namely $a_0>0$ (where we are assuming nondegeneracy). Therefore, in the
Routh-Hurwitz case three equations are obtained. The first two are:
\begin{eqnarray}
\Re(\lambda_v)< 0\,,\\
\Re(\lambda_v) \left[ |\lambda_v|^2-2\Re(\lambda_x)\right]
-\left[\Re(\lambda_x)\Re(\lambda_v)+\Im(\lambda_x)\Im(\lambda_v)\right] > 0\,.
\end{eqnarray}
The third inequality we do not utilize, since it is extremely complicated
containing fifth order terms. We are left with the above two, which are now
necessary conditions for all roots to have negative real part.

Similarly, the Li\'{e}nard-Chipart stability criterion also gives two
necessary conditions for all roots to have negative real part:
\begin{eqnarray}
\Re(\lambda_v) < 0\,,\\
2\Re(\lambda_x)-|\lambda_v|^2 < 0\,.
\end{eqnarray}
The third inequality is the same as before and will not be utilized. Since the
second inequality of the Li\'{e}nard-Chipart conditions seems less complicated
than the corresponding one of the Routh-Hurwitz conditions, we will continue
with the former.

Substituting the expressions for the $\lambda$'s (Lemma~\ref{lemm:the lambdas})
we get:
\begin{eqnarray*}
g_v \left[\sum_{j=0}^{2} \alpha_{v,j}\cos(j\phi)\right] &<& 0\,, \\
g_x \left[\sum_{j=0}^{2} \alpha_{x,j}\cos(j\phi)\right]-
g_v^2\left\{\left[\sum_{j=0}^{2} \alpha_{v,j}\cos(j\phi)\right]^2-
\left[\sum_{j=0}^{2}\beta_{v,j}\sin(j\phi)\right]^2\right\}&<& 0\,.
\end{eqnarray*}
These are complicated relations therefore we will use the equivalent relations
averaged over $\phi$. The first of these equations was already used in
Lemma~\ref{lemm:simple conditions}. After some calculations we can work out the
average over $\phi$ of the second relation. This gives the final necessary
condition for all roots to have negative real part.

\begin{cory}
The following is a necessary condition for the collection ${\cal S}^*$ to not be b
unstable:
\begin{eqnarray*}
g_{x}-g_{v}^2\sum_{j=-2}^{2} \rho_{v,j}^2\le 0\,.
\end{eqnarray*}
\label{cor:rhc}
\end{cory}

\section*{Appendix C: Analysis of Type~II Trajectories}

\begin{prop} Let $K_0>0$ fixed. Suppose that ${\cal S}^*$ is stable and that
$0<c_-<c_+$. Suppose further that there are $\alpha\in(0,1)$ and $q>0$ such
that $Na_m m^{1+q}$ and $Nb_m m^{1+q}$ are bounded, and $(2-q)\alpha\leq 1$.
Then
\begin{eqnarray*}
\lim_{N\rightarrow\infty} \sup_{t\in[0,K_0N]}\,|z_N(t)-{\overline z}_N(t)|=0\,.
\end{eqnarray*}
where ${\overline z}_N(t)$ is given by
\begin{eqnarray*}
{\overline z}_N(t)= \left\{\begin{array}{cl}-t & t\in \left[0,\frac{N}{c_+}\right)\\
\frac{-N}{c_+} + \left[\frac{c_-}{c_+ - c_-}\right)\left(t-\frac{N}{c_+}\right) &
t\in\left(\frac{N}{c_+},\frac{N}{c_-}\right)\\
0 & t\in\left[\frac{N}{c_-}, \infty\right)
\end{array} \right.
\end{eqnarray*}
The signal velocities are as in Theorem~\ref{thm:velocities}.
\end{prop}
{\bf Sketch of Proof:} We consider the equations of motion for the
acceleration $\xi_k$ of agent $k$. These are given by the second derivative
with respect to time of Definition~\ref{def:nnn}. In those equations, the only
expression that depends on time is the initial condition of leader. So nothing
changes, except that now $\xi_0(t)=\delta(t)$ (for $>0$), where $\delta$ is
the Dirac function. We replace the Dirac function by a smooth pulse $p(t)$
that enables us to satisfy the decay constraint on the decay of $a_m$ and
$b_m$ but with the condition that $\int\,p(s)\,ds=1$. So now we obtain:
\begin{eqnarray}
\xi_0(t) = p(t)
\label{eq:pulse}
\end{eqnarray}
Theorem~\ref{thm:velocities} now implies that in ${\cal S}^*$ we have
\begin{eqnarray}
\xi_k(t)= f_+(k-c_+t)+f_-(k-c_-t)
\label{eq:soln-acceleration}
\end{eqnarray}
By the periodic boundary conditions conjectures, we see that away from the
boundaries the behavior of ${\cal S}$ and ${\cal S}^*$ is the same. So we have
the above relation from $t=0$ until the signal runs into the boundary at $N$.

Setting $k=0$ in the last equation and comparing with
Equation~(\ref{eq:pulse}) gives
\begin{eqnarray}
p(t) = f_+(-c_+t)+f_-(-c_-t)\,.
\label{eq:substitution}
\end{eqnarray}
The second part of Equation~(\ref{eq:newBC}) then gives:
\begin{eqnarray*}
f_+'(-c_+t)+f_-'(-c_-t)=0\,.
\end{eqnarray*}
Integrate with respect to $t$ to get
\begin{eqnarray*}
\frac{-1}{c_+}f_+(-c_+t)-\frac{1}{c_-}f_-(-c_-t)=0 \quad \Longrightarrow \quad
f_-(s)=-\frac{c_-}{c_+}f_+\left(\frac{c_+}{c_-}\,s\right)\,.
\end{eqnarray*}
Substitute this into Equation~(\ref{eq:substitution}):
\begin{eqnarray*}
p(t)=\frac{c_+-c_-}{c_+}f_+(-c_+t)\quad \Longrightarrow \quad
f_+(s)=\frac{c_+}{c_+-c_-}p\left(\frac{s}{c_+}\right)\,.
\end{eqnarray*}
Now use both of the last equations to eliminate $f_-$ and $f_+$ from
Equation~(\ref{eq:soln-acceleration}):
\begin{eqnarray*}
\xi_k(t)=\frac{c_+}{c_+-c_-}\, p\left(t-\frac{k}{c_+}\right) -
\frac{c_-}{c_+-c_-}\, p\left(t-\frac{k}{c_-}\right)\,.
\end{eqnarray*}

Now set $k=N$ and integrate twice with respect to time and add a Galilean
transformation so that for small positive $t$ we get $z_N(t)=-t$. With some
rewriting this gives the final result. (As before, to actually prove the
remainder indeed tends to zero, one needs the assumption on the decay of the
$a_m$ and $b_m$. This part of the argument is given in \cite{Cantos2014-1}.)
\qed

\bibliographystyle{unsrt}
\bibliography{manunnn}

\end{document}